\documentclass[10pt]{article}
\usepackage{amsmath}
\usepackage[all,matrix,arrow,curve,poly,arc]{xy}

\newcommand{\qed}{%
  \ifmmode 
   \eqno{\qedsymbol}
  \else
    \leavevmode\unskip\penalty9999 \hbox{}\nobreak\hfill\hbox{\qedsymbol}
  \fi
}
\newcommand{\qedsymbol}{\leavevmode\vrule height 1.2ex width 1.1ex depth -.1ex}
\newenvironment{proof}{\begin{trivlist}\item[\hskip%
\labelsep{\bf Proof.\quad}]}%
{\hfill\qed\rm\end{trivlist}}

\newtheorem{theorem}{Theorem}[section]
\newtheorem{corollary}[theorem]{Corollary}
\newtheorem{lemma}[theorem]{Lemma}
\newtheorem{proposition}[theorem]{Proposition}

\def\r{\;{\cal R}\;}
\def\l{\;{\cal L}\;}
\def\rs{\;{\cal R^\ast}\;}
\def\ls{\;{\cal L^\ast}\;}

\def\x{{\overline{x}}}
\def\y{{\overline{y}}}
\def\exy{{\overline{xy}}}

\def\e{\overline{e}}

\begin{document}
\date{}
\title{Adequate transversals of quasi-adequate semigroups}
\author{{\Large Jehan Al-Bar and James Renshaw\footnote{Communicating author}}\\School of Mathematics\\University of Southampton\\Southampton, SO17 1BJ\\England\\Email: j.h.renshaw@maths.soton.ac.uk\\jaal\_bar@hotmail.com}
\maketitle

\begin{abstract}
The concept of an adequate transversal of an abundant semigroup was introduced by El-Qallali in~\cite{el-qallali} whilst in \cite{fountain3}, he and Fountain initiated the study of quasi-adequate semigroups as natural generalisations of orthodox semigroups. In this work we provide a structure theorem for adequate transversals of certain types of quasi-adequate semigroup and from this deduce Saito's classic result on the structure of inverse transversals of orthodox semigroups. We also apply it to left ample adequate transversals of left adequate semigroups and provide a structure for these based on semidirect products of adequate semigroups by left regular bands.
\end{abstract}

\medskip


\medskip


\section{Introduction and preliminaries}\label{introduction-section}

The concept of an inverse transversal was essentially introduced in 1978 by Blyth and McAlister~\cite{mcalister-blyth} in response to investigations into the structure of split orthodox semigroups and then later by Blyth and McFadden to natural questions concerning the greatest idempotent in a naturally ordered regular semigroup. Let $S^0$ be an inverse subsemigroup of a regular semigroup $S$ and suppose that for all $x\in S, |S^0\cap V(x)| = 1$. Then $S^0$ is called an {\em inverse transversal} of $S$ and we denote the unique inverse of $x$ in $S^0$ by $x^0$.
Over the years a great deal of work has been done to examine the structure of regular semigroups containing inverse transversals and to consider a number of useful and important properties that these transversals may possess. The reader is referred to the excellent survey article in~\cite{blyth} for more details. As one might expect a number of generalisations of this concept have emerged and connections and similarlities with inverse transversals studied. For example, we may replace the inverse subsemigroup with a different type of regular subsemigroup such as an orthodox semigroup and, with a suitable adjustment of the definitions, refer to {\em orthodox transversals} of regular semigroups (see for example~\cite{chen2}). Or we may replace the unique inverse $x^0$ of the regular element $x$ with another related element such as an {\em associate}  - recall that an {\em associate} of an element $x$ is an element $y$ such that $xyx = x$ - and study the structure of these {\em associate inverse} subsemigroups (see for example~\cite{billhardt}).
In this work we are more concerned with what was probably the first of these generalisations, initiated by El-Qallali, to {\em adequate transversals} of {\em abundant} semigroups and in particular to a structure theorem for quasi-adequate semigroups with an adequate transversal.

\smallskip

The focus of attention of this paper is on structure theorems for adequate transversals mainly following the standard ``Rees Theorem'' type approach. We do however also consider a structure based on spined products which has already been shown to be useful in the inverse transversal case. A number of authors, such as Chen (\cite{chen}), have studied the structure of certain types of adequate transversal but a general structure theorem seem to be rather too difficult to contemplate. Consequently we shall consider only the structure of certain types of quasi-adequate semigroup and in some cases restrict our attention to certain types of transversal, such as quasi-ideal transversals.

\medskip

After some preliminary results and concepts in section 1 we concentrate on quasi-adequate semigroups in section 2 and in particular on those quasi-adequate semigroups $S$ with a splitting morphism $S\to S/\gamma$ where $\gamma$ is the least admissible adequate congruence on $S$. Quasi-adequate semigroups are analoguous to orthodox semigroups in the regular case. We culminate in this section with a rather technical structure theorem for adequate transversals of such semigroups which is of a similar nature to one given by Saito for orthodox semigroups. From this we deduce a number of very interesting corollaries concerning adequate transversals which are quasi-ideals one of which involves a spined product of two important subsemigroups of the quasi-adequate semigroup first introduced in~\cite{albar-renshaw}. In section 3 we focus our attention on the structure of left ample adequate transversals of quasi-adequate semigroups which are also left adequate and describe such semigroups in terms of a rather elegant semidirect-product construction. Finally in section 4 we consider what happens when our abundant semigroups are in fact regular and deduce a number of new as well as known structure theorems for inverse transversals.

\medskip

Unless otherwise stated the terminology and notation will be that of~\cite{howie}. Let $S$ be a semigroup and define a left congruence ${\cal R}^\ast$ on $S$ by
$$
\rs=\{(a,b)\in S\times S\; |\; xa=ya \text{ if and only if } xb=yb \text{ for all } x,y\in S^1\}.
$$
The right congruence ${\cal L}^\ast$ is defined dually and we shall let ${\cal H}^\ast = {\cal L}^\ast\cap{\cal R}^\ast$. It is easy to show that if $a$ and $b$ are regular elements of $S$ then $a\rs b$ if and only if $a\r b$. We say that a semigroup is {\em abundant} if each ${\cal R}^\ast-$class and each ${\cal L}^\ast-$class contains an idempotent. An abundant semigroup in which the idempotents commute is called {\em adequate}. It is clear that regular semigroups are abundant and that inverse semigroups are adequate.

\begin{lemma}{\rm(\cite[Corollary 1.2]{fountain2})}\label{e-rs-lemma}
Let $e\in E(S)$ and $a \in S$. Then $e\rs a$ if and only if $ea=a$ and for all $x,y\in S^1$, $xa=ya$ implies $xe=ye$.
\end{lemma}

\begin{lemma}{\rm(\cite[Proposition 1.3]{fountain1})}
A semigroup $S$ is adequate if and only if each ${\cal L}^\ast-$class and each ${\cal R}^\ast-$class contain a unique idempotent and the subsemigroup generated by $E(S)$ is regular.
\end{lemma}

If $S$ is an adequate semigroup and $a\in S$ then we shall denote by $a^\ast$ the unique idempotent in $L_a^\ast$ and by $a^+$ the unique idempotent in $R_a^\ast$. It is easy to show that if $S$ is adequate and $a,b\in S$ then $a\rs b$ if and only if $a^+=b^+$ and $a\ls b$ if and only if $a^\ast=b^\ast$.

\begin{lemma}{\rm(\cite[Proposition 1.6]{fountain1})}\label{star-plus-lemma}
If $S$ is an adequate semigroup then for all $a,b\in S,\ (ab)^\ast = (a^\ast b)^\ast$ and $(ab)^+=(ab^+)^+$.
\end{lemma}

Notice that we can then immediately deduce that for all $a,b \in S, a^+(ab)^+ = (ab)^+$ and that $(ab)^\ast b^\ast = (ab)^\ast$.

\medskip

If $S$ is an abundant semigroup and $U$ is an abundant subsemigroup of $S$ then we say that $U$ is a {\em $\ast-$subsemigroup} of $S$ if ${\cal L}^\ast(U) = {\cal L}^\ast(S)\cap(U\times U), {\cal R}^\ast(U) = {\cal R}^\ast(S)\cap(U\times U)$. It can be shown that $U$ is a $\ast-$subsemigroup of $S$ if and only if for all $a \in U$ there exist $e,f \in E(U)$ such that $e\in L^\ast_a(S), f\in R^\ast_a(S)$ (see~\cite{el-qallali}).

\bigskip

Now suppose that $S^0$ is an adequate $\ast-$subsemigroup of the abundant semigroup $S$. We say that $S^0$ is an {\em adequate transversal} of $S$ if for each $x \in S$ there is a unique $\x\in S^0$ and $e,f\in E$ such that
$$
x = e\x f \text{ and such that } e\l\x^+ \text{ and } f\r\x^\ast.
$$
It is straightforward to show,~\cite{el-qallali}, that such an $e$ and $f$ are uniquely determined by $x$. Hence we normally denote $e$ by $e_x$, $f$ by $f_x$ and the semilattice of idempotents of $S^0$ by $E^0$.

\begin{lemma}{\rm(\cite[Lemma 1.4]{albar-renshaw})}\label{e-f-lemma}
Let $S$ be an abundant semigroup with an adequate transversal $S^0$. Then for all $x\in S$
\begin{enumerate}
\item $e_x\rs x$ and $f_x\ls x$,
\item if $x \in S^0$ then $e_x = x^+\in E^0, \x=x, f_x = x^\ast\in E^0$,
\item if $x \in E^0$ then $e_x = \x = f_x = x$,
\item $e_\x\l e_x$ and hence $e_\x e_x = e_\x$ and $e_x e_\x = e_x$,
\item $f_\x\r f_x$ and hence $f_\x f_x = f_x$ and $f_x f_\x = f_\x$.
\end{enumerate}
\end{lemma}

\medskip

\begin{lemma}\ {\rm\cite[Proposition 2.3]{chen}}\label{rs-ls-lemma}
Let $S^0$ be an adequate transversal of an abundant semigroup $S$ and let $x,y\in S$. Then
\begin{enumerate}
\item $x\rs y$ if and only if $e_x=e_y$,
\item $x\ls y$ if and only if $f_x=f_y$.
\end{enumerate}
\end{lemma}

\bigskip

We define
$$
I = \{e_x : x \in S\},\quad\Lambda=\{f_x : x \in S\}
$$
and note from the previous results that there are bijections $I\to S/{\cal R^\ast}$ and $\Lambda\to S/{\cal L^\ast}$. It is also straightforward to see that if $x \in I, y\in \Lambda$ then
$$e_x = x, \x=f_x=e_\x\text{ and }e_y=\y=f_\y, f_y=y.$$
Also, it is well-known that $|R^\ast_x\cap I| = 1$ and that $|L^\ast_x\cap\Lambda|=1$, results that will prove useful later.

\medskip

Suppose now that $x\in Reg(S)$, the set of regular elements of $S$. Using the fact that $x\r e_x$ and $x\l f_x$ then from \cite[Theorem 2.3.4]{howie} there exists a unique $x^0\in V(x)$ with $xx^0=e_x$ and $x^0x=f_x$.
In what follows we shall write $x^{00}$ for $\left(x^0\right)^0$.

\begin{theorem}\ {\rm\cite[Theorems 2.3 \& 2.4]{albar-renshaw}}\label{main-regular-theorem}
Let $S^0$ be an adequate transversal of an abundant semigroup. If $x\in Reg(S)$ then $|V(x)\cap S^0|=1$. Moreover $x^0\in S^0, \x = x^{00}$ and $x^0 = x^{000}$.
Also,
$$I = \{x\in Reg(S) : x = xx^0\} = \{xx^0 : x\in Reg(S)\}$$
and
$$\Lambda = \{x\in Reg(S) : x = x^0x\} = \{x^0x : x \in Reg(S)\}.$$
\end{theorem}

\medskip

In addition, from Lemma~\ref{rs-ls-lemma} we can easily deduce

\begin{lemma}
Let $S^0$ be an adequate transversal of an abundant semigroup. For all $x \in I, a \in Reg(S)$, $x = aa^0$ if and only if $a\in {\cal R}^\ast_x\cap Reg(S)$.
\end{lemma}
Consequently we see that
\begin{lemma}
Let $S^0$ be an adequate transversal of an abundant semigroup and let $x\in S, a\in Reg(S)$. Then $e_x = aa^0$ and $f_x = a^0a$ if and only if $a\in H_x^\ast=R_x^\ast\cap L_x^\ast$.
\end{lemma}

\medskip

Notice that $R_x^\ast\cap Reg(S)\ne\emptyset$ and $L_x^\ast\cap Reg(S)\ne\emptyset$ as $e_x\in R_x^\ast\cap Reg(S)$ and $f_x\in L_x^\ast\cap Reg(S)$.

\medskip

Let $T = Reg(S)$, let $U = T\cap S^0$ and suppose that $T$ is a subsemigroup of $S$. It is clear that $U$ is then a regular subsemigroup of $T$. It is also relatively straightforward to see that $U = \{x^0 : x \in T\} = \{x\in T : x = \x = x^{00}\}$. Given that $S^0$ is adequate then we can easily deduce that $U$ is an inverse transversal of the regular semigroup $T$.
Moreover, it is reasonably clear that, with the obvious notation, $E(T) = E(S), E(U) = E(S^0)$ and $I(T) = I(S), \Lambda(T) = \Lambda(S)$ and so from~\cite[Theorem 1.3]{blyth} we see that

\begin{proposition}\ {\rm\cite[Proposition 2.5 \& Proposition 2.6]{albar-renshaw}}\label{subband-proposition}
Let $S^0$ be an adequate transversal of an abundant semigroup. If $T$ is a subsemigroup of $S$ then $I$ (resp. $\Lambda$) is a left (resp. right) regular subband of $S$
and for all $x,y \in T$
$$
(xy)^0 = (x^0xy)^0x^0 = y^0(xyy^0)^0 = y^0(x^0xyy^0)^0x^0,
$$
and
$$
(xy^0)^0 = y^{00}x^0, (x^0y)^0 = y^0x^{00}.
$$
\end{proposition}

\section{Quasi-Adequate Semigroups}

A semigroup is said to be {\em quasi-adequate} if it is abundant and its idempotents form a subsemigroup. It was shown in ~\cite[Proposition 1.3]{fountain3} that in this case the set $T$ of regular elements is an orthodox subsemigroup of $S$.
From section~\ref{introduction-section} we see that $U = T\cap S^0$ is an inverse transversal of $T$. The following is now immediate from~\cite[Theorem 1.12]{blyth}.

\begin{proposition}{\rm(\cite[Proposition 6.2]{albar-renshaw})}\label{quasi-adequate-proposition}
Let $S^0$ be an adequate transversal of an abundant semigroup $S$. Then the following are equivalent:
\begin{enumerate}
\item $S$ is quasi-adequate;
\item $(\forall x,y\in T)$, $(xy)^0 = y^0x^0$;
\item $(\forall i\in I)(\forall l\in\Lambda)$, $(li)^0 = i^0l^0$;
\item $I\Lambda = E(S)$.
\end{enumerate}
\end{proposition}

\smallskip

\begin{corollary}\label{e0-corollary}
Let $S^0$ be an adequate transversal of a quasi-adequate semigroup $S$. If $x \in E(S)$ then $x^0\in E^0$. In particular if $x\in I\cup\Lambda$ then $x^0\in E^0$.
\end{corollary}

\begin{lemma}\label{E0-transversal-lemma}
Let $S^0$ be an adequate transversal of a quasi-adequate semigroup $S$. Then $I$ (resp. $\Lambda$) is a left (resp. right) regular subband of $S$ and $E^0$ is a semilattice transversal of both $I$ and $\Lambda$ with $x^0\in V(x)\cap E^0$ for all $x\in I\cup\Lambda$. Moreover if $x\in E^0$ then $x^0 = x$.
\end{lemma}

\begin{proof}
That $I$ is a left regular subband of $S$ follows from Proposition~\ref{subband-proposition}. 
We prove that $E^0$ is a transversal of $I$, the result for $\Lambda$ being similar. For $x\in I$ we know that $x^0\in V(x)\cap E^0$. Suppose that $x'\in V(x)\cap E^0$. Then $xx'\in I$ since $I$ is a band. Hence $xx'\rs x$ and so $xx' = e_x = xx^0$ since $|R^\ast_x\cap I|= 1$. Similarly, $x'x = x^0x$ and so by the definition of $x^0$ it follows that $x' = x^0$ as required.
\end{proof}

\smallskip

The following two subsets of $S$ are defined in~\cite{albar-renshaw} and the diagram below summaries the various connections between these sets.
$$
R = \{x \in S : e_x = e_\x\},\ L = \{x \in S : f_x = f_\x\}
$$
$$
\xy 0;/r10pc/:
{\ellipse<30pt>{}};p;p+(0.25,0),{\ellipse<30pt>{}};
p;p+(-0.125,-.25),{\ellipse<30pt>{}};
p+(0,0.25)*+!D{S^0},+(-0.4,0.05)*+!RD{L},+(0.8,0)*+!LD{R},+(-0.4,-0.55)*+!U{E(S)},+(-.12,0.35)*+!R{I},,+(.24,0)*+!L{\Lambda}
\endxy
$$
From~\cite[Theorem 3.1]{albar-renshaw} we see that $R = \{x \in S: x = \x f_x\}, L = \{x \in S: x = e_x\x\}$.
\begin{proposition}{\rm(\cite[Theorem 3.11 \& Corollary 3.13]{albar-renshaw})}\label{bar-proposition}
Let $S$ be an abundant semigroup with an adequate transversal $S^0$ and suppose that $x,y \in S$. Then $\exy = \x\;\y$ in each of the following situations:
\begin{enumerate}
\item $x\in \Lambda, y \in I$;
\item $x\in L, y \in R$;
\item $x,y \in S^0$.
\end{enumerate}
\end{proposition}

From Theorem~\ref{main-regular-theorem} and Proposition~\ref{quasi-adequate-proposition} we can also deduce

\begin{proposition}
Let $S$ be a quasi-adequate semigroup with an adequate transversal $S^0$ and suppose that $x,y \in T=Reg(S)$. Then $\exy = \x\;\y$.
\end{proposition}

\begin{corollary}
Let $S$ be an orthodox semigroup with an adequate (and hence inverse) transversal $S^0$. Then for all $x,y\in S$, $\exy = \x\;\y$.
\end{corollary}

We say that $S^0$ is a {\em quasi-ideal} of $S$ if $S^0SS^0\subseteq S^0$ or equivalently~\cite[Proposition 2.2]{chen} if $\Lambda I\subseteq S^0$ or~\cite[Lemma 1.4]{albar-renshaw2} if $RL\subseteq S^0$. These transversals have been the subject of a great deal of study in both the inverse and adequate cases.

\begin{proposition}
Let $S$ be an abundant semigroup with a quasi-ideal adequate transversal $S^0$. $S$ is quasi-adequate if and only if for all $x,y \in S, \exy = \x\;\y$.
\end{proposition}

\begin{proof}
First note that if $S^0$ is a quasi-ideal then from~\cite[Theorem 3.12]{albar-renshaw} we see that for all $x,y\in S$
$$\exy = \x f_xe_y\y,\quad e_{xy} = e_xe_{xy},\quad f_{xy} = f_{xy}f_x.$$

\smallskip

Suppose that $S$ is quasi-adequate.
From Proposition~\ref{quasi-adequate-proposition} we deduce that $\overline{f_xe_y} = f_\x e_\y$ and since $S^0$ is a quasi-ideal of $S$ then $f_xe_y = \overline{f_xe_y}$. Hence $\exy = \x f_xe_y\y = \x f_\x e_\y\y = \x\;\y$.

\medskip
Conversely suppose that $\exy=\x\;\y$ for all $x,y\in S$. Let $x,y\in E$ so that $\x,\y\in E^0$ and hence $\x\;\y \in E^0$. Since $S^0$ is a quasi-ideal of $S$ then $\y\;\x = \overline{yx} = \y f_ye_x\x$. Similarly $\x\;\y = \x f_xe_y\y$. Hence
$$
xyxy = xe_y\y f_ye_x\x f_xy = xe_y\y\;\x f_xy = e_x\x f_x e_y\y\;\x f_x e_y\y f_y = e_x\x\;\y\;\x\;\y f_y
$$
$$= e_x\x\;\y f_y = e_x\x f_xe_y\y f_y = xy.
$$
So $S$ is quasi-adequate as required.
\end{proof}

\medskip

Let $S$ be a quasi-adequate semigroup with band of idempotents $E$ and for $e\in E$ let $E(e)$ denote the ${\cal J}-$class of $e$ in $E$. For $a\in S$, let $a^+$ denote a typical element of $R^\ast_a(S)\cap E$ and let $a^\ast$ denote a typical element of $L^\ast_a(S)\cap E$. Define a relation $\delta$ on $S$ by
$$
\delta=\{(a,b)\in S\times S : b = eaf, \text{ for some } e\in E(a^+), f\in E(a^\ast)\}.
$$
From~\cite{fountain3} we see that $\delta$ is an equivalence relation and is contained in any adequate congruence $\rho$ on $S$. Recall~\cite{fountain3} that a morphism $\phi : S\to T$ is called {\em good} if for all $a,b \in S$, $a\;{\cal R}^\ast(S)\;b$ implies $a\phi\;{\cal R}^\ast(T)\;b\phi$ and $a\;{\cal L}^\ast(S)\;b$ implies $a\phi\;{\cal L}^\ast(T)\;b\phi$. So for example we see that if $U$ is an abundant subsemigroup of an abundant semigroup $S$ then $U$ is a $\ast-$subsemigroup of $S$ if and only if the inclusion $U\to S$ is good. Recently the term {\em admissible} has been used in place of {\em good} and we shall henceforth use that term. A congruence $\rho$ is called {\em admissible} if the natural homomorphism $\rho^\natural : S\to S/\rho$ is admissible. It was shown in~\cite[Proposition 2.1]{fountain3} that if $S$ is quasi-adequate then there exists a minimum adequate admissible congruence $\gamma$ on $S$.

\begin{lemma}{\ \rm(\cite[Proposition 2.6]{fountain3})}
If $\delta$ is a congruence then $\delta$ is the minimum adequate admissible congruence on $S$.
\end{lemma}

An {\em idempotent-connected} (IC) abundant semigroup $S$, is an abundant semigroup in which for each $a \in S$ and some $a^+\in R^\ast_a\cap E, a^\ast\in L^\ast_a\cap E$, there is a bijection $\alpha:\langle a^+\rangle\to\langle a^\ast\rangle$ such that $xa = a(a\alpha)$ for all $a\in\langle a^+\rangle$. A {\em bountiful} semigroup is an IC quasi-adequate semigroup.

\begin{theorem}{\ \rm(\cite[Theorem 2.6]{guo})}
Let $S$ be a bountiful semigroup. Then $\delta$ is an admissible congruence.
\end{theorem}

\begin{theorem}{\ \rm(\cite[Corollary 2.8]{fountain3})}
Let $S$ be a quasi-adequate semigroup with band of idempotents $E$. If $E$ is normal then $\delta$ is an admissible congruence.
\end{theorem}

The following characterisation of the property that $\exy=\x\;\y$ for all $x,y\in S$ essentially appears in~\cite{el-qallali}.

\begin{proposition}{\ \rm(\cite[Proposition 4.3 \& Proposition 4.4]{el-qallali})}
If $S$ is a quasi-adequate semigroup with an adequate transversal $S^0$ then the following are equivalent
\begin{enumerate}
\item $\delta$ is a congruence on $S$,
\item $\delta=\{(a,b)\in S\times S: \overline{a}=\overline{b}\}$,
\item for all $x,y \in S$, $\exy = \x\;\y$.
\end{enumerate}
 Moreover in this case $S/\delta\cong S^0$.
\end{proposition}
\medskip

Consequently, we shall say that an adequate transversal $S^0$ of a quasi-adequate semigroup $S$ is {\em admissible} if $\exy = \x\;\y$ for all $x,y \in S$. Notice that for such a transversal the natural map $\delta^\natural:S\to S^0$ {\em splits} in the sense that $\iota\delta^\natural = 1_{S^0}$, where $\iota:S^0\to S$ is the inclusion morphism. Hence $S$ is a split quasi-adequate semigroup.

\begin{lemma}\label{xybar-lemma}
Let $S$ be a quasi-adequate semigroup with an admissible adequate transversal $S^0$. If $x,y \in S$
then $\exy = \overline{\x f_xe_y\y}$.
\end{lemma}

\begin{proof}
Since $\exy = \x\;\y$ then we see that
$$\exy = \overline{xf_x}\;\overline{e_yy} = \x\;\overline{f_x}\;\overline{e_y}\;\y = \overline{\x}\;\overline{f_x}\;\overline{e_y}\;\overline{\y} = \overline{\x f_xe_y\y}.$$
\end{proof}

\begin{theorem}\label{quasi-adequate-ef-theorem}
Let $S$ be a quasi-adequate semigroup with an admissible adequate transversal $S^0$. Then $e_{xy} = e_xe_{\x f_xe_y\y}$ and $f_{xy} = f_{\x f_xe_y\y}f_y$.
\end{theorem}

\begin{proof}
Notice first that using Lemma~\ref{xybar-lemma}
\begin{align*}
\left(e_xe_{\x f_xe_y\y}\right)\left(\exy\right)\left(f_{\x f_xe_y\y}f_y\right)
&= e_xe_{\x f_xe_y\y}\overline{\x f_xe_y\y}f_{\x f_xe_y\y}f_y\\
&= e_x\x f_xe_y\y f_y\\
&=xy
\end{align*}
Moreover $e_x\l\x^+$ and $e_{\x f_xe_y\y}\l \overline{\x f_xe_y\y}^+ = \left(\exy\right)^+$ from Lemma~\ref{xybar-lemma} and so we have
\begin{align*}
\left(e_xe_{\x f_x e_y\y}\right)(\exy)^+&= e_x\left(e_{\x f_x e_y\y}(\exy)^+\right)\\
& = e_xe_{\x f_x e_y\y}
\end{align*}
and from Lemmas~\ref{e-rs-lemma},~\ref{star-plus-lemma} and~\ref{e-f-lemma}
\begin{align*}
(\exy)^+\left(e_xe_{\x f_x e_y\y}\right)& = (\x^+\exy)^+\left(e_xe_{\x f_x e_y\y}\right)\\
& = (\x^+(\exy)^+)^+\left(e_xe_{\x f_x e_y\y}\right)\\
& = (\exy)^+\x^+\left(e_xe_{\x f_x e_y\y}\right)\\
& = (\exy)^+\x^+\left(e_{\x f_x e_y\y}\right)\\
& =  (\exy)^+e_{\x f_x e_y\y}\\
& =  (\exy)^+.
\end{align*}
So $e_xe_{\x f_x e_y\y}\l(\exy)^+$.
In a similar way $f_{\x f_xe_y\y}f_y\r(\exy)^\ast$ and the result follows.
\end{proof}

We therefore see that if $S$ is a quasi-adequate semigroup with an admissible adequate transversal $S^0$ then we have the factorisation
$$
(e_x\x f_x)(e_y\y f_y) = (e_xe_{\x f_xe_y\y})\left(\overline{\x f_xe_y\y}\right)(f_{\x f_xe_y\y}f_y).
$$

Let $I$ be a left regular band with semilattice transversal $E^0$ and let $x \in I$. Denote the ${\cal L}-$class of $x$ by $L_x$ and notice that $L_x = \{y \in I: y^0 = x^0\} = L_{x^0}$. To see this notice that if $x\l y$ then $x^0x = y^0y$ and so $x^0 = x^0xx^0 = x^0x = y^0y = y^0yy^0 = y^0$. Also $x^{00} = x^0$ since $E^0$ is a semilattice and so $x^0\l x$.

\begin{lemma}
For all $x \in E^0$, $L_x$ is a left zero semigroup and if $x,y \in E^0$ then $L_yL_x\subseteq L_{yx}$.
\end{lemma}

\begin{proof}
Let $x \in E^0$. Then $L_x$ is left zero since if $a,b\in L_x$ then $a\l b$ and since $a,b$ are idempotents then $ab=a$.
If $x,y\in E^0$ and if $c\in L_y, d\in L_x$ then $(cd)^0 = d^0c^0 = x^0y^0 = (yx)^0 = x^0$ and so $L_yL_x\subseteq L_{yx}$.
\end{proof}

Hence $I = \cup_{x\in E^0}{L_x}$ and in this case we see that $I$ is a semilattice $E^0$ of the left zero semigroups $L_x$.
In a similar way, for each $x \in \Lambda$,  $R_x$ is a right zero semigroup and if $x\le y$ then $R_xR_y\subseteq R_x$. Hence $\Lambda = \cup_{x\in E^0}{R_x}$ and we see that $\Lambda$ is a semilattice $E^0$ of the right zero semigroups $R_x$. Notice that the above result is not true in general as can be seen from~\cite[Example 2.7]{chen}. 

\medskip

It is worth noting that by Lemma~\ref{star-plus-lemma}, for all $x,y \in S^0$ we have $(xy)^\ast \le x^\ast$ and so $R_{(xy)^\ast}R_{y^\ast}\subseteq R_{(xy)^\ast}$ and $(xy)^+\le x^+$ so that $L_{x^+}L_{(xy)^+}\subseteq L_{(xy)^+}$.

\medskip

\begin{lemma}\label{bar-x+-lemma}
Let $S$ be an abundant semigroup with an adequate transversal $S^0$. Then for all $x\in S^0, a\in \Lambda, b\in I$, it follows that $a\in R_{x^\ast}$ if and only if $\overline{a} = x^\ast$ whilst $b\in L_{x^+}$ if and only if $\overline{b} = x^+$.
\end{lemma}
\begin{proof}
First notice that if $a \in \Lambda, x \in S^0$ and if $a \in R_{x^\ast}$ then $x^\ast\r a\r \overline{a}\in E^0$ and so $x^\ast = \overline{a}$ as $E^0$ is sub-semilattice of $S$. Conversely if $x^\ast = \overline{a}$ then $a\r\overline{a} = x^\ast$ and so $a \in R_{x^\ast}$. The result for $b\in I$ is similar.
\end{proof}

\medskip

We now come to the main result of the paper.
\begin{theorem}\label{structure-theorem}
Let $S^0$ be an adequate semigroup with semilattice $E^0$ and let $I = \cup_{x\in E^0}{L_x}$ be a left regular band and $\Lambda = \cup_{x\in E^0}{R_x}$ a right regular band with a common semilattice transversal $E^0$. Suppose that for each $x,y \in S^0$ there exist mappings $\alpha_{x,y} : R_{x^\ast}\times L_{y^+}\to L_{(xy)^+}
$ and $\beta_{x,y} : R_{x^\ast}\times L_{y^+}\to R_{(xy)^\ast}
$ satisfying:
\begin{enumerate}
\item if $f \in R_{x^\ast}, g \in L_{y^+}, h \in R_{y^\ast}, k \in L_{z^+}$ then
$$(f,g)\alpha_{x,y}\left((f,g)\beta_{x,y}h,k\right)\alpha_{xy,z} = \left(f,g(h,k)\alpha_{y,z}\right)\alpha_{x,yz}$$
$$\left(f,g(h,k)\alpha_{y,z}\right)\beta_{x,yz}(h,k)\beta_{y,z} = \left((f,g)\beta_{x,y}h,k\right)\beta_{xy,z},$$
\item $(x^\ast,y^+)\alpha_{x,y} = (xy)^+, (x^\ast,y^+)\beta_{x,y} = (xy)^\ast$,
\item if $x, x_1, x_2\in S^0, e_1\in L_{x_1^+}, f_1 \in R_{x_1^\ast}, e_2\in L_{x_2^+}, f_2 \in R_{x_2^\ast}, e\in L_{x^+}$ and if
$$e_1(f_1,e)\alpha_{x_1,x} = e_2(f_2,e)\alpha_{x_2,x},\ x_1x = x_2x \text{ and } (f_1,e)\beta_{x_1,x}x^\ast = (f_2,e)\beta_{x_2,x}x^\ast$$
then
$$e_1(f_1,e)\alpha_{x_1,x^+} = e_2(f_2,e)\alpha_{x_2,x^+},\ x_1x^+ = x_2x^+ \text{ and } (f_1,e)\beta_{x_1,x^+} = (f_2,e)\beta_{x_2,x^+}.$$
\item if $x,x_1, x_2\in S^0, e_1\in L_{x_1^+}, f_1 \in R_{x_1^\ast}, e_2\in L_{x_2^+}, f_2 \in R_{x_2^\ast}, f\in R_{x^\ast}$ and if
$$x^+(f,e_1)\alpha_{x,x_1} = x^+(f,e_2)\alpha_{x,x_2},\ xx_1 = xx_2 \text{ and } (f,e_1)\beta_{x,x_1}f_1 = (f,e_2)\beta_{x,x_2}f_2$$
then
$$(f,e_1)\alpha_{x^\ast,x_1} = (f,e_2)\alpha_{x^\ast,x_2},\ x^\ast x_1 = x^\ast x_2 \text{ and } (f,e_1)\beta_{x^\ast,x_1}f_1 = (f,e_2)\beta_{x^\ast,x_2}f_2$$
\end{enumerate}
Define a multiplication on the set
$$
W=\{(e,x,f)\in I\times S^0\times\Lambda : e \in L_{x^+}, f \in R_{x^\ast}\}
$$
by
$$(e,x,f)(g,y,h) = (e(f,g)\alpha_{x,y},xy,(f,g)\beta_{x,y}h).
$$
Then $W$ is a quasi-adequate semigroup with an admissible adequate transversal isomorphic to $S^0$
. Moreover, if in addition $\alpha$ and $\beta$ satisfy:
\begin{enumerate}
\item[5.] for all $f\in\Lambda, e \in I$,
$$(f^0,e)\alpha_{f^0,e^0} = f^0e,\quad (f,e^0)\beta_{f^0,e^0} = fe^0,$$
\end{enumerate}
then $I(W) \cong I, \Lambda(W)\cong \Lambda$.

Moreover every quasi-adequate semigroup $S$, with an admissible adequate transversal can be constructed in this way.
\end{theorem}

\begin{proof}
We establish the result in stages. First we prove that $W$ is a quasi-adequate semigroup. Next we show that $W$ has an adequate transversal isomorphic to $S^0$. Then we demonstrate that this transversal is admissible and finally we consider what happens when property (5) holds.

\medskip

That $W$ is a semigroup follows easily from property 1. Notice now that $(e,x,f) \in E(W)$ if and only if $x\in E^0$. To see this, suppose that $(e,x,f)\in E(W)$. Then
$$
(e,x,f) = (e,x,f)(e,x,f) = (e(f,e)\alpha_{x,x},x^2,(f,e)\beta_{x,x}f),
$$
and so $x = x^2$. Conversely, if $x = x^2$ then
$$
(e,x,f)(e,x,f) = (e(f,e)\alpha_{x,x},x^2,(f,e)\beta_{x,x}f) = (e,x,f)
$$
since $(f,e)\alpha_{x,x} \in L_{{(x^2)}^+} = L_{x^+}$ and so since $L_{x^+}$ is a left zero semigroup then $e(f,e)\alpha_{x,x} = e$. Similarly $(f,e)\beta_{x,x}f=f$.

Now we show that $(e,x^+,x^+)\rs(e,x,f)$. First notice that
$$
(e,x^+,x^+)(e,x,f) = (e(x^+,e)\alpha_{x^+,x},x^+x,(x^+,e)\beta_{x^+,x}f) = (e,x,f).
$$
Suppose then that $(e_1,x_1,f_1)(e,x,f) = (e_2,x_2,f_2)(e,x,f)$. Then
$$
(e_1(f_1,e)\alpha_{x_1,x},x_1x,(f_1,e)\beta_{x_1,x}f) = (e_2(f_2,e)\alpha_{x_2,x},x_2x,(f_2,e)\beta_{x_2,x}f).
$$
Since $f\r x^\ast$ then $fx^\ast = x^\ast$ and so $(f_1,e)\beta_{x_1,x}x^\ast = (f_2,e)\beta_{x_2,x}x^\ast$. Hence from (3) we see that $e_1(f_1,e)\alpha_{x_1,x^+} = e_2(f_2,e)\alpha_{x_2,x^+}$, 
$x_1x^+ = x_2x^+$ and $(f_1,e)\beta_{x_1,x^+} = (f_2,e)\beta_{x_2,x^+}$ which gives
$$
\begin{array}{rl}
(e_1,x_1,f_1)(e,x^+,x^+)& = (e_1(f_1,e)\alpha_{x_1,x^+},x_1x^+,(f_1,e)\beta_{x_1,x^+}x^+)\\& = (e_2(f_2,e)\alpha_{x_2,x^+},x_2x^+,(f_2,e)\beta_{x_2,x^+}x^+)\\& = (e_2,x_2,f_2)(e,x^+,x^+)\\
\end{array}
$$
as required. In a similar way and using (4), $(e,x,f)\ls(x^\ast,x^\ast,f)$. Consequently we see that $W$ is an abundant semigroup. In addition, if $(e,x,f),(g,y,h) \in E(W)$ then $x,y \in E^0$ and $(e,x,f)(g,y,h) = (e(f,g)\alpha_{x,y},xy,(f,g)\beta_{x,y}h)$. Hence $(e,x,f)(g,y,h) \in E(W)$ and $W$ is quasi-adequate.

\medskip

Now we show that $W$ has an adequate transversal isomorphic to $S^0$.

\medskip

To this end let $W^0 = \{(x^+,x,x^\ast) : x \in S^0\}$. First notice that $W^0$ is a subsemigroup of $W$ since by (3) we have
$$
(x^+,x,x^\ast)(y^+,y,y^\ast) = (x^+(x^\ast,y^+)\alpha_{x,y},xy,(x^\ast,y^+)\beta_{x,y}y^\ast) = ((xy)^+,xy,(xy)^\ast).
$$
This also clearly demonstrates that $W^0$ is isomorphic to $S^0$ and therefore is adequate. Now it is clear that $E(W^0) = \{(x,x,x) : x \in E^0\}$ and so from the above arguments, we see that $(x^+,x,x^\ast)^+ = (x^+,x^+,x^+)$ and that $(x^+,x,x^\ast)^\ast = (x^\ast,x^\ast,x^\ast)$. It also follows from the arguments above that $(x^+,x^+,x^+){\cal R}^\ast_W(x^+,x,x^\ast)$ and that $(x^\ast,x^\ast,x^\ast){\cal L}^\ast_W(x^+,x,x^\ast)$ and so $W^0$ is a $\ast-$subsemigroup of $W$. To see that $W^0$ is an adequate transversal of $W$ notice that
$$
(e,x,f) = (e,x^+,x^+)(x^+,x,x^\ast)(x^\ast,x^\ast,f)
$$
and from above we have that $(e,x^+,x^+)\l(x^+,x,x^\ast)^+, (x^\ast, x^\ast,f)\r(x^+,x,x^\ast)^\ast$. We need only demonstrate that these terms are unique with respect to these properties. So suppose that
$$
(e,x,f) = (g,z,h)(y^+,y,y^\ast)(j,w,k)
$$
with $(g,z,h),(j,w,k)\in E(W)$ and $(g,z,h)\l(y^+,y,y^\ast)^+=(y^+,y^+,y^+),$ $(j,w,k)\r(y^+,y,y^\ast)^\ast=(y^\ast,y^\ast,y^\ast)$. Then $(g,z,h)(y^+,y^+,y^+) = (g,z,h)$ and $(y^+,y^+,y^+)(g,z,h) = (y^+,y^+,y^+)$ and so $z\l y^+$. In a similar manner $w\r y^\ast$ and so since $S^0$ is adequate we see that $z=y^+$ and $w=y^\ast$. Hence $x = zyw = y$ as required. In fact it is not too hard to deduce that $g=e, k=f, h=x^+$ and $j=x^\ast$.

\medskip

To see that $W^0$ is an admissible transversal of $W$ notice that if $(e,x,f),(g,y,h) \in W$ then it is straightforward to check, using (2), that
$$\overline{(e,x,f)(g,y,h)} = ((xy)^+,xy,(xy)^\ast) = \overline{(e,x,f)}\;\overline{(g,y,h)}.$$

\medskip

Suppose now that for all $f\in\Lambda, e \in I$, $(f^0,e)\alpha_{f^0,e^0} = f^0e$, and $(f,e^0)\beta_{f^0,e^0} = fe^0$. It is worth noting that if $(e,x,f)\in W$ then from Corollary~\ref{e0-corollary} and Lemma~\ref{bar-x+-lemma} we see that $e^0 = e^{00} = \overline{e} = x^+$ and so with the obvious notation
$$
I(W) = \{(e,e^0,e^0) : e\in I\}.
$$
In a similar way
$$
\Lambda(W) = \{(f^0,f^0,f) : f \in \Lambda\}.
$$
So for $e,g\in I$
$$
(e,e^0,e^0)(g,g^0,g^0) = (e(e^0,g)\alpha_{e^0,g^0},e^0g^0,(e^0,g)\beta_{e^0,g^0}g^0)
$$
Since we know that $W^0$ is an adequate transversal of $W$ and that $W$ is quasi-adequate then it follows that $I(W)$ is a band. Hence $(e^0,g)\beta_{e^0,g^0}g^0 = e^0g^0$ and so we see
$$
\begin{array}{rl}
(e,e^0,e^0)(g,g^0,g^0)& = (e(e^0,g)\alpha_{e^0,g^0},e^0g^0,(e^0,g)\beta_{e^0,g^0}g^0)\\
&= (ee^0g,e^0g^0,e^0g^0)\\
&= (eg,(eg)^0,(eg)^0).
\end{array}
$$
Hence $I(W)\cong I$. Similarly $\Lambda(W)\cong\Lambda$.
\medskip

Conversely, let $S$ be a quasi-adequate semigroup with an admissible adequate transversal $S^0$. Let $I$ and $\Lambda$ be as in the adequate transversal decomposition and note that by Lemma~\ref{E0-transversal-lemma}, $I$ is a left regular subband and $\Lambda$ is a right regular subband of $S$ with a common semilattice transversal $E^0$.
Define $W = \{(e_x,\x,f_x) : x\in S\}$ and let $S \to W$ be given by $x\mapsto(e_x,\x,f_x)$. Then since $e_x\l\x^+$ and $f_x\l\x^\ast$ it follows that $W$ is of the required form. For $x,y \in S^0$ and $a\in R_{x^\ast}, b \in L_{y^+}$ define $(a,b)\alpha_{x,y} = e_{xaby}$ and $(a,b)\beta_{x,y} = f_{xaby}$ so that by Theorem~\ref{quasi-adequate-ef-theorem} the multiplication in $W$ is given by
$$(e_x,\x,f_x)(e_y,\y,f_y) = (e_xe_{\x f_xe_y\y},\x\;\y,f_{\x f_xe_y\y}f_y) = (e_{xy},\exy,f_{xy})$$
which is clearly associative and so $W$ is a semigroup and the map $S\to W$ is a morphism.
It is also clearly a bijection and so an isomorphism.
We show that $\alpha$ and $\beta$ satisfy the properties given in the statement of the theorem.

First notice that $(a,b)\alpha_{x,y} = e_{xaby}\l\overline{(xaby)}^+ = \left(\overline{x}\;\overline{a}\;\overline{b}\;\overline{y}\right)^+=\left(x\overline{a}\;\overline{b}y\right)^+ = (xx^\ast y^+y)^+ = (xy)^+$ as required. Similarly, $(a,b)\beta_{x,y}\r(xy)^\ast$.

\medskip

\begin{enumerate}
\item Let $a \in R_{x^\ast}, b \in L_{y+}, c \in R_{y^\ast}, d \in L_{z^+}$. Then
$$(a,b)\alpha_{x,y}\left((a,b)\beta_{x,y}c,d)\right)\alpha_{xy,z} = e_{xaby}\left(f_{xaby}c,d\right)\alpha_{xy,z}=e_{xaby}e_{xyf_{xaby}cdz}$$
and
$$\left(a,b(c,d)\alpha_{y,z}\right)\alpha_{x,yz}=\left(a,be_{ycdz}\right)\alpha_{x,yz}=e_{xabe_{ycdz}yz}.$$
Now if we put $u = xaby, v = cdz$  then we see that $\overline{u} = \overline{xaby} = \x\;\overline{a}\;\overline{b}\;\y = xx^\ast y^+y = xy$. In a similar way, $\overline{ycdz} = yz$. From Theorem~\ref{quasi-adequate-ef-theorem} it follows that
$$
\begin{array}{rl}
e_{uv}& = e_{(uv)f_v}\\
&= e_{uv}e_{\overline{uv}f_{uv}e_{f_v}f_{\overline{v}}}\\
&= e_ue_{\overline{u}f_ue_v\overline{v}}e_{\overline{uv}f_{\overline{u}f_ue_v\overline{v}}f_ve_{f_v}f_{\overline{v}}}\\
&= e_u\left(e_{\overline{u}f_ue_v\overline{v}}e_{\overline{\overline{u}f_ue_v\overline{v}}f_{\overline{u}f_ue_v\overline{v}}e_{f_v}f_{\overline{v}}}\right)\\
&= e_ue_{(\overline{u}f_ue_v\overline{v})f_v}\\
&= e_ue_{\overline{u}f_uv}\\
&= e_{xaby}e_{xyf_{xaby}cdz}.
\end{array}
$$
On the other hand we notice that
$\overline{e_{ycdz}yz} = \overline{yz} = yz$ and $e_{(e_{ycdz}yz)} = e_{(ycdz)}e_{yz}$. Then
$$
\begin{array}{rl}
e_{((xab)e_{(ycdz)}yz)}&= e_{xab}e_{(\overline{xab}f_{(xab)}e_{(e_{ycdz}yz)}\overline{yz})}\\
&= e_{xab}e_{(\overline{xab})f_{(xab)}e_{(ycdz)}e_{(yz)}\overline{ycdz}}\\
&= e_{xab}e_{(\overline{xab})f_{(xab)}e_{(ycdz)}e_{(\overline{yz})}\overline{ycdz}}\\
&= e_{xab}e_{(\overline{xab})f_{(xab)}e_{(ycdz)}e_{(\overline{ycdz})}\overline{ycdz}}\\
&= e_{xab}e_{(\overline{xab})f_{(xab)}e_{(ycdz)}\overline{ycdz}}\\
&= e_{(xab)(ycdz)}\\
&= e_{uv}.
\end{array}
$$

\medskip
In a similar manner we can establish that 
$$\left(f,g(h,k)\alpha_{y,z}\right)\beta_{x,yz}(h,k)\beta_{y,z} = \left((f,g)\beta_{x,y}h,k\right)\beta_{xy,z}.$$
\item $(x^\ast,y^+)\alpha_{x,y} = e_{xx^\ast y^+y} = e_{xy} = (xy)^+, (x^\ast,y^+)\beta_{x,y} = f_{xx^\ast y^+y} = f_{xy} = (xy)^\ast$,
\item First, since $x_1x=x_2x$ and since $x\rs x^+$ then $x_1x^+=x_2x^+$. Also, $(f_1,e)\alpha_{x_1,x} = e_{x_1f_1ex} \rs x_1f_1ex \rs x_1fex^+ \rs e_{x_1f_1ex^+} = (f_1,e)\alpha_{x_1,x^+}$ and so $(f_1,e)\alpha_{x_1,x} = (f_1,e)\alpha_{x_1,x^+}$ (since $|R^\ast_{x_1f_1ex}\cap I|=1$) and so
$$e_1(f_1,e)\alpha_{x_1,x^+}=e_1(f_1,e)\alpha_{x_1,x}=e_2(f_2,e)\alpha_{x_2,x}=e_2(f_2,e)\alpha_{x_2,x^+}.$$
Finally, since $(f_1,e)\beta_{x_1,x}x^\ast = (f_2,e)\beta_{x_2,x}x^\ast$ then $f_{x_1f_1ex}x^\ast =f_{x_2f_2ex}x^\ast$. Now $f_{x_1f_1ex}\ls x_1f_1ex$ and so $f_{x_1f_1ex}x^\ast\ls x_1f_1exx^\ast = x_1f_1ex\ls f_{x_1f_1ex}$. Consequently $f_{x_1f_1ex}x^\ast = f_{x_1f_1ex}$ since $x^\ast\in E^0\subseteq\Lambda$ ($|\Lambda\cap L^\ast_{x_1f_1ex}| = 1$). It follows that $f_{x_1f_1ex} = f_{x_2f_2ex}$. Hence we see that $e_1e_{x_1f_1ex}x_1xf_{x_1f_1ex} = e_2e_{x_2f_2ex}x_2xf_{x_2f_2ex}$ and so from Theorem~\ref{quasi-adequate-ef-theorem} we see that $e_1x_1f_1ex = e_2x_2f_2ex$. Since $x\rs x^+$ then $e_1x_1f_1ex^+ = e_2x_2f_2x^+$. Now let $y = e_1x_1f_1$ and $z = ex^+x^+$ and note from Theorem~\ref{quasi-adequate-ef-theorem} that $f_{yz} = f_{\overline{y}f_ye_z\overline{z}}f_z = f_{x_1f_1ex^+}x^+$. In a similar way, if $w = e_2x_2f_2$ then $f_{wz} = f_{\overline{w}f_we_z\overline{z}}f_z = f_{x_2f_2ex^+}x^+$. Hence we deduce that $f_{x_1f_1ex^+}x^+ = f_{x_2f_2ex^+}x^+$. Finally $f_{x_1f_1ex^+}\ls x_1f_1ex^+$ and so $f_{x_1f_1ex^+}x^+\ls x_1f_1ex^+\ls f_{x_1f_1ex^+}$. Therefore $f_{x_1f_1ex^+}\ls f_{x_2f_2ex^+}$ and so $(f_1,e)\beta_{x_1,x^+} = f_{x_1f_1ex^+} = f_{x_2f_2ex^+} = (f_2,e)\beta_{x_2,x^+}$ as required.
\item Similar to (3).
\item By definition $(f^0,e)\alpha_{f^0,e^0} = e_{f^0f^0ee^0} = e_{f^0e} = f^0e$ since $f^0e\in E^0I\subseteq I$. In a similar way $(f,e^0)\beta_{f^0,e^0} = fe^0$.
\end{enumerate}
\end{proof}

\medskip

As an example of Theorem~\ref{structure-theorem}, suppose that $S^0$ is an adequate semigroup with semilattice of idempotents $E^0$ and let $I = \cup_{x\in E^0}{L_x}$ be a left normal band and $\Lambda = \cup_{x\in E^0}{R_x}$ a right normal band with a common quasi-ideal semilattice transversal $E^0$. For each $x,y \in S^0$ define $\alpha_{x,y} : R_{x^\ast}\times L_{y^+} \to L_{(xy)^+}$ and $\beta_{x,y} : R_{x^\ast}\times L_{y^+}\to R_{(xy)^\ast}$ by $(f,g)\alpha_{x,y} = (xy)^+, (f,g)\beta_{x,y} = (xy)^\ast$. Then
\begin{enumerate}
\item if $f \in R_{x^\ast}, g \in L_{y^+}, h \in R_{y^\ast}, k \in L_{z^+}$ then
$$(f,g)\alpha_{x,y}\left((f,g)\beta_{x,y}h,k\right)\alpha_{xy,z} = (xy)^+((xy)z)^+ = (x(yz))^+ = \left(f,g(h,k)\alpha_{y,z}\right)\alpha_{x,yz}$$
and
$$\left(f,g(h,k)\alpha_{y,z}\right)\beta_{x,yz}(h,k)\beta_{y,z} = (x(yz))^\ast(yz)^\ast = (x(yz))^\ast = \left((f,g)\beta_{x,y}h,k\right)\beta_{xy,z},$$
\item $(x^\ast,y^+)\alpha_{x,y} = (xy)^+, (x^\ast,y^+)\beta_{x,y} = (xy)^\ast$,
\item if $x, x_1, x_2\in S^0, e_1\in L_{x_1^+}, f_1 \in R_{x_1^\ast}, e_2\in L_{x_2^+}, f_2 \in R_{x_2^\ast}, e\in L_{x^+}$ and if
$$e_1(f_1,e)\alpha_{x_1,x} = e_2(f_2,e)\alpha_{x_2,x},\ x_1x = x_2x \text{ and } (f_1,e)\beta_{x_1,x}x^\ast = (f_2,e)\beta_{x_2,x}x^\ast$$
then
$$e_1(x_1x)^+ = e_2(x_2x)^+,\ x_1x = x_2x \text{ and } (x_1x)^\ast x^\ast = (x_2x)^\ast x^\ast$$
and so since $x\rs x^+$, $(x_1x^+)^+ = (x_1x)^+$ and since $(x_1x)^\ast x^\ast = (x_1x)^\ast$ then
$$e_1(x_1x^+)^+ = e_2(x_2x^+)^+,\ x_1x^+ = x_2x^+ \text{ and } (x_1x^+)^\ast = (x_2x^+)^\ast$$
and hence
$$e_1(f_1,e)\alpha_{x_1,x^+} = e_2(f_2,e)\alpha_{x_2,x^+},\ x_1x^+ = x_2x^+ \text{ and } (f_1,e)\beta_{x_1,x^+} = (f_2,e)\beta_{x_2,x^+}.$$
\item in a similar way, if $x,x_1, x_2\in S^0, e_1\in L_{x_1^+}, f_1 \in R_{x_1^\ast}, e_2\in L_{x_2^+}, f_2 \in R_{x_2^\ast}, f\in R_{x^\ast}$ and if
$$x^+(f,e_1)\alpha_{x,x_1} = x^+(f,e_2)\alpha_{x,x_2},\ xx_1 = xx_2 \text{ and } (f,e_1)\beta_{x,x_1}f_1 = (f,e_2)\beta_{x,x_2}f_2$$
then
$$(f,e_1)\alpha_{x^\ast,x_1} = (f,e_2)\alpha_{x^\ast,x_2},\ x^\ast x_1 = x^\ast x_2 \text{ and } (f,e_1)\beta_{x^\ast,x_1}f_1 = (f,e_2)\beta_{x^\ast,x_2}f_2$$
\item if $f\in\Lambda, e \in I$ then since $E^0$ is a quasi-ideal of $I$ we have
$$(f^0,e)\alpha_{f^0,e^0} = (f^0e^0)^+ = f^0e^0 = f^{000}e^{00} = (f^0e)^{00} = \overline{f^0e} = f^0e
$$
and in a similar way
$$
(f,e^0)\beta_{f^0,e^0} = (f^0e^0)^\ast = f^0e^0 = f^{00}e^{000} = (fe^0)^{00} = \overline{fe^0} = fe^0. 
$$
\end{enumerate}
Consequently we see from Theorem~\ref{structure-theorem} that the set
$$
W=\{(e,x,f)\in I\times S^0\times\Lambda : e \in L_{x^+}, f \in R_{x^\ast}\}
$$
together with multiplication defined by
$$(e,x,f)(g,y,h) = (e(xy)^+,xy,(xy)^\ast h).
$$
is a quasi-adequate semigroup with an admissible adequate transversal $W^0=\{(x^+,x,x^\ast) : x \in S^0\}$ isomorphic to $S^0$. In addition, since property (5) holds, we see that $I(W) = \{(e,x^+,x^+) : x \in S^0, e\in L_{x^+}\}\cong I$ and $\Lambda(W) = \{(x^\ast,x^\ast,f) : x \in S^0, f \in R_{x^\ast}\}\cong \Lambda$ and so
$$
(x^\ast,x^\ast,f)(e,y^+,y^+) = (x^\ast(x^\ast y^+)^+,x^\ast y^+,(x^\ast y^+)^\ast y^+) = (x^\ast y^+,x^\ast y^+,x^\ast y^+)\in W^0
$$
which means that $W^0$ is a quasi-ideal of $W$.

\medskip

On the other hand, if $S$ is a quasi-adequate semigroup with an admissible adequate transversal $S^0$ such that $S^0$ is a quasi-ideal of $S$ then we see from Theorem~\ref{structure-theorem} that for all $x,y \in S^0$ there exists maps $\alpha_{x,y}$ and $\beta_{x,y}$ with the properties (1)-(5) and such that
\begin{center}
$S\cong W = \{(e,x,f)\in I\times S^0\times\Lambda : e \in L_{x^+}, f \in R_{x^\ast}\}$ and $S^0\cong W^0 = \{(x^+,x,x^\ast) : x \in S^0\}$
\end{center}
with multiplication given by
$$
(e,x,f)(g,y,h) = (e(f,g)\alpha_{x,y},xy,(f,g)\beta_{x,y} h).
$$
Moreover from the proof we see that $(f,g)\alpha_{x,y} = e_{xfgy} = e_{\overline{xfgy}} = (\overline{xfgy})^+ = (\x\overline{f}\overline{g}\;\y)^+ = (\x\;\y)^+ = (xy)^+$. In a similar way $(f,g)\beta_{x,y} = (xy)^\ast$. Also, if $f\in E^0, e \in I$ then by property (5), $fe = f^0e = (f^0,e)\alpha_{f^0,e^0} = (f^0e^0)^+ \in E^0$ and so $E^0$ is a quasi-ideal of $I$. In a similar way $E^0$ is a quasi-ideal of $\Lambda$. Finally, we know from \cite[Proposition 1.7]{mcalister} together with the remarks before Proposition~\ref{subband-proposition} that $I$ is a left normal band and $\Lambda$ is a right normal band and so we have established

\begin{corollary}\label{quasi-ideal-corollary}
Let $S^0$ be an adequate semigroup with semilattice of idempotents $E^0$ and let $I = \cup_{x\in E^0}{L_x}$ be a left normal band and $\Lambda = \cup_{x\in E^0}{R_x}$ be a right normal band with a common semilattice transversal $E^0$. Let
$$
W=\{(e,x,f)\in I\times S^0\times\Lambda : e \in L_{x^+}, f \in R_{x^\ast}\}
$$
and define a multiplication on $W$ by
$$
(e,x,f)(g,y,h) = (e(xy)^+,xy,(xy)^\ast h).
$$
Then $W$ is a quasi-adequate semigroup with a quasi-ideal, admissible adequate transversal isomorphic to $S^0$. Conversely every such transversal can be constructed in this way.
\end{corollary}

An adequate transversal $S^0$ of an abundant semigroup $S$ is said to be {\em multiplicative} if $\Lambda I\subseteq E(S^0)$. It is worth noting that by~\cite[Theorem 5.3 \& Corollary 6.3]{albar-renshaw} that if $S$ is quasi-adequate then $S^0$ is multiplicative if and only if $S^0$ is a quasi-ideal. Notice also that since $I$ is left normal then by~\cite[Lemma 1.7~(3) \& Theorem 4.2]{albar-renshaw} $E^0$ is a quasi-ideal of both $I$ and $\Lambda$.

\bigskip

An interesting consequence of this corollary is the following alternative characterisation involving spined products.
An abundant semigroup is said to be {\em left (resp. right) adequate} if every ${\cal R}^\ast-$class (resp. ${\cal L}^\ast-$class) contains a unique idempotent. From~\cite[Theorem 3.14]{albar-renshaw} we see that is $S$ is a left adequate semigroup with an adequate transversal $S^0$ then $\Lambda = E^0, R = S^0, L = S$ and $I = E(S)$.
It also follows from~\cite[Theorem 2.3]{albar-renshaw2} that if $S^0$ is a quasi-ideal adequate transversal of an abundant semigroup $S$ then $L$ and $R$ are subsemigroups of $S$, $L$ is left adequate and $R$ is right adequate with $S^0$ a common quasi-ideal adequate transversal of both $L$ and $R$.

\def\a{\overline{a}}
\def\b{\overline{b}}
\def\c{\overline{c}}
\def\d{\overline{d}}
\def\h{\overline{h}}
\def\l{\overline{l}}
\def\spd{L|\times|R}

\medskip

Suppose now that $S^0$ is a quasi-ideal adequate transversal of an abundant semigroup $S$.
Notice that by~\cite[Corollary 3.13]{albar-renshaw}, if $a \in L, b\in R, x \in S^0$ then
\begin{center}
$\overline{ax} = \a\;\x, e_{ax} = e_a(\a\;\x)^+, f_{ax}=(\a\;\x)^\ast f_x$ and $\overline{xb} = \x\;\b, e_{xb} = e_x(\x\;\b)^+, f_{xb} = (\x\;\b)^\ast f_b$.
\end{center}
In particular, if $a\in L_{x^+}, b\in R_{x^\ast}$ then by Lemma~\ref{bar-x+-lemma}, $\a = x^+, \b = x^\ast$ and so
\begin{center}
$\overline{ax} = x, e_{ax} = a, f_{ax} = f_x$ and $\overline{xb}=x, e_{xb} = e_x, f_{xb} =b$.
\end{center}

Since $S^0$ is a quasi-ideal of $S$ then $RL\subseteq S^0$ and so we can define a multiplication on the spined product
$$
\spd = \{(x,a)\in L\times R : \x = \a\}
$$
by
$$
(x,a)(y,b) = (x\y,\a b) = (x\b,\x b)
$$
and it is an easy matter to demonstrate that under this multiplication $\spd$ is a semigroup.
\begin{corollary}\label{spined-product-corollary}
Let $L$ be a left adequate semigroup and $R$ a right adequate semigroup with a common quasi-ideal adequate transversal $S^0$. Construct the spined product
$$
\spd = \{(x,a)\in L\times R : \x = \a\}
$$
and define a multiplication on $\spd$ by
$$
(x,a)(y,b) = (x\y,\a b) = (x\b,\x b)
$$
Then $\spd$ is a quasi-adequate semigroup with an admissible, quasi-ideal adequate transversal isomorphic to $S^0$. Moreover every such transversal can be constructed in this way.
\end{corollary}

\begin{proof}
Let $L,R$ and $S^0$ be as in the statement of the theorem, let $E^0$ be the semilattice of idempotents of $S^0$ and let $I= E(L)=I(L)$ and $\Lambda = E(R) = \Lambda(R)$.
To show that $I$ is left normal suppose that $e\in E^0, i \in I(L) = E(L)$. Then $ei \in E^0I = \Lambda(L) I \subseteq S^0$ since $S^0$ is a quasi-ideal of $L$ and so $ei = \overline{ei} = \overline{e}\;\overline{i} = e\;\overline{i}$ by Proposition~\ref{bar-proposition}~(1). Hence by~\cite[Theorem 4.2]{albar-renshaw} we see that $I$ is left normal. In a similar way $\Lambda$ is right normal.
Since $S^0$ is a common adequate transversal of $L$ and $R$ then it follows easily that $E^0$ is a common semilattice transversal of $I$ and $\Lambda$.
Now from Corollary~\ref{quasi-ideal-corollary}
$$
W=\{(g,x,l)\in I\times S^0\times\Lambda : g \in L_{x^+}, l \in R_{x^\ast}\}
$$
is a quasi-adequate semigroup with multiplication given by
$$
(g,x,l)(h,y,k) = (g(xy)^+,xy,(xy)^\ast k)
$$
and $W^0 = \{(x^+,x,x^\ast) x \in S^0\} \cong S^0$ is a quasi-ideal, admissible adequate transversal of $W$. Consider then the map $\phi : W \to \spd$ given by $(g,x,l) \mapsto (gx,xl)$. Then
$$
\begin{array}{ll}
(g,x,l)(h,y,k)\phi&=(g(xy)^+,xy,(xy)^\ast k)\phi = (g(xy),(xy)k)\\
&=(g(x\overline{hy}),(\overline{xl}y)k)=(gx,xl)(hy,yk)\\
&=(g,x,l)\phi(h,y,k)\phi
\end{array}
$$
and so $\phi$ is a morphism. If $x\in L, a \in R$ then $(x,a) = (e_x\x,\a f_a) = (e_x,\x,f_a)\phi$ and hence $\phi$ is onto while if $(g,x,l)\phi = (h,y,k)\phi$ then $(gx,xl) = (hy,yk)$ and so from above we see that $x = \overline{gx} = \overline{hy} = y$. In addition, $g = e_{gx} = e_{hy} = h$ and $l = f_{xl} = f_{yk} = k$ and so $\phi$ is an isomorphism. Clearly $W^0 = \{(x^+,x,x^\ast) x \in S^0\}\cong S^0$.

\medskip

Conversely suppose that $S$ is a quasi-adequate semigroup with a quasi-ideal, admissible adequate transversal $S^0$ of $S$ and suppose that $E^0$ is the semilattice of idempotents of $S^0$. Then by Corollary~\ref{quasi-ideal-corollary} there exists a left normal band $I$ and a right normal band $\Lambda$ with a common semilattice transversal $E^0$ such that
$$
W=\{(g,x,l)\in I\times S^0\times\Lambda : g \in L_{x^+}, l \in R_{x^\ast}\}
$$
is a quasi-adequate semigroup with multiplication given by
$$
(g,x,l)(h,y,k) = (g(xy)^+,xy,(xy)^\ast k)
$$
and in addition $I(W) \cong I=I(S), \Lambda(W)\cong \Lambda=\Lambda(S)$. Infact we see from the proof that
$$
W=\{(e_x,\x,f_x)\in I\times S^0\times\Lambda : x\in S\}.
$$
Let $L = L(S), R = R(S)$ and consider the map $\theta : \spd \to W$ given by $(x,a)\theta = (e_x,\x,f_a)$. Then
$(x,a)(y,b)\theta = (x\y,\a b)\theta = (e_{x\y},\overline{x\y},f_{\a b})=(e_x(\x\;\y)^+,\x\;\y,(\a\;\b)^\ast f_b)=(e_x,\x,f_a)(e_y,\y,f_b)=(x,a)\theta(y,b)\theta$ and so $\theta$ is a morphism which is clearly onto. If $(x,a)\theta = (y,b)\theta$ then $(e_x,\x,f_a) = (e_y,\y,f_b)$ and so $e_x = e_y, \x=\y$ and $f_x = f_\x = f_{\a} = f_{\b} = f_\y = f_y$ and so $x = y$.
In a similar way $a = b$ and so $\theta$ is an isomorphism. That $L$ is left adequate and $R$ is right adequate follow from the statements above.
\end{proof}

The interested reader may like to consult~\cite{albar-renshaw2} for results of a similar nature involving spined product decompositions of abundant semigroups.

\smallskip

In view of~\cite[Section 7]{albar-renshaw} it would be of interest to know what effect $S$ being a monoid has on the details of Theorem~\ref{structure-theorem}.

\section{Left Adequate Semigroups}

As another application of Theorem~\ref{structure-theorem} we consider the case of left adequate semigroups. Let $S$ be a left adequate semigroup with an adequate transversal $S^0$ and as usual let $E^0$ be the semilattice of idempotents of $S^0$. Since $\Lambda = E^0$ then if $S$ is also quasi-adequate and $S^0$ is admissible then we must have $R_{x^\ast} = \{x^\ast\}$ and so $(f,e)\beta_{x,y} = (xy)^\ast$. In addition we can also deduce from \cite[Theorem 3.14]{albar-renshaw} that for all $x \in S, x = e_x\x$ and we shall make use of this fact in what follows without further reference.

\bigskip
Let $S^0$ be an adequate semigroup with semilattice of idempotents $E^0$ and suppose that $I$ is a left regular band with a semilattice transversal (isomorphic to) $E^0$. Suppose also that there is defined on $I$ a left $S^0-$action $S^0\times I \to I$ given by $(x,e) \mapsto x\ast e$ and which is distributive over the multiplication on $I$. In other words, for all $x,y \in S^0, e,f \in I$ we have
$$
(xy)\ast e = x\ast(y\ast e)\text{ and\ }x\ast(ef) = (x\ast e)(x\ast f).
$$
We can construct the semidirect product of $S^0$ by $I$ as 
$$
I\ast S^0 = \{(e,x)\in I\times S^0\}
$$
with multiplication given by
$$
(e,x)(g,y) = (e(x\ast g),xy)
$$
and it is an easy matter to check that $I\ast S^0$ is a semigroup.
Now consider the subsemigroup $W = \{(e,x) \in I\ast S^0: e \in L_{x^+}\}$.

\begin{lemma}\label{ast-lemma}
Let $I,S^0$ and $W$ be as above and suppose that for all $x,y \in S^0, x\ast y^+ = (xy)^+$. Then
\begin{enumerate}
\item for all $x,y \in S^0, f \in L_{y^+}$, $x\ast f\in L_{(xy)^+}$;
\item $E(W) = \{(e,x)\in W : x\in E^0\}$;
\item $E(W)$ is a band.
\end{enumerate}
\end{lemma}

\begin{proof} Let $x,y \in S^0$ and $e \in I$.
\begin{enumerate}
\item Let $f\in L_{y^+}$ and notice that
$$
\begin{array}{rcl}
x\ast f &=&x\ast(fy^+) = (x\ast f)(x\ast y^+)= (x\ast f)(xy)^+,\\
(xy)^+(x\ast f)&=&(x\ast y^+)(x\ast f) = x\ast(y^+f) = x\ast y^+ = (xy)^+,\\
\end{array}
$$
and so $x\ast f \in L_{(xy)^+}$ as required.

\item If $(e,x) \in E(W)$ then $(e,x) = (e,x)(e,x) = (e(x\ast e),x^2)$ and so $x \in E^0$. Conversely, if $x\in E^0$ then from the first part, $x\ast e \in L_{(x^2)^+} = L_{x^+}$ and so $e\l x\ast e$ and $e = e(x\ast e)$. Consequently $(e,x) \in E(W)$.
\item That $E(W)$ is a band is reasonably clear.
\end{enumerate}
\end{proof}

\begin{lemma}\label{left-adequate-lemma}
Let $I,S^0$ and $W$ be as above and suppose that for all $x,y \in S^0, x\ast y^+ = (xy)^+$. Then $W$ is a left abundant semigroup and for all $w\in W, |R^\ast_{w}\cap E(W)| = 1$.
\end{lemma}

\begin{proof}
Let $(e,x) \in W$ and notice that $(e,x^+)(e,x) = (e(x^+\ast e),x^+x) = (e,x)$ since $x^+\ast e \in L_{x^+}$ by the previous Lemma. In addition, if $(e_1,x_1)(e,x) = (e_2,x_2)(e,x)$ then $(e_1(x_1\ast e),x_1x) = (e_2(x_2\ast e),x_2x)$ and since $x\rs x^+$ then $(e_1,x_1)(e,x^+) = (e_2,x_2)(e,x^+)$ and so $(e,x^+)\rs(e,x)$ and $W$ is left abundant.

Suppose now that $(f,y) \in E(W)\cap R^\ast_{(e,x)}$. Then $(f,y)(e,x) = (e,x)$ and so $(f(y\ast e),yx) = (e,x)$. Hence $yx = x$ and $f(y\ast e) = e$. It follows since $x\rs x^+$ that $yx^+ = x^+$. Now it is easy to check, using Lemma~\ref{ast-lemma}(1), that $(e,x^+)(e,x) = (f,y)(e,x)$ and so it follows that $(e,x^+)(f,y) = (f,y)(f,y)$ and hence $x^+y = y$. Consequently we deduce that $x^+ = y$. But from the Lemma~\ref{ast-lemma}(1) $y\ast e \in L_{(yx)^+} = L_{x^+} = L_{y^+}$ and so $e=f(y\ast e) = f$ as required.
\end{proof}

\medskip

Let $W^0 = \{(x^+,x) : x \in S^0\}$ and notice that $W^0$ is a subsemigroup of $W$. From the proof of the previous Lemma we see that $(x^+,x^+)\rs(x^+,x)$ and so by~\cite[Proposition 1.3]{el-qallali}, $W^0$ is a right $\ast-$subsemigroup of $W$ and hence is also  left abundant.

\medskip

We say that a left adequate semigroup $S$ is {\em left ample} (formerly called {\em left type-A}) if for all $a \in S, e \in E(S), (ae) = (ae)^+a$.

\begin{lemma}\label{left-ample-lemma}
Let $I,S^0$ and $W$ be as above and suppose that for all $x,y \in S^0, x\ast y^+ = (xy)^+$ and that $S^0$ is left ample. Then for all $w\in W^0, e \in E(W^0), we = (we)^+w$.
\end{lemma}
\begin{proof}
Let $w = (x^+,x)$ and let $e = (y^+,y^+)$ for some $x,y \in S^0$. Then $we = (x^+,x)(y^+,y^+) = ((xy)^+,xy^+)$ and $(we)^+ = ((xy)^+,(xy)^+)$. Hence $(we)^+w = ((xy)^+,(xy)^+)(x^+,x) = ((xy)^+,(xy)^+x) = we$ as required.
\end{proof}

\begin{theorem}\label{left-adequate-theorem}
Let $S^0$ be a left ample, adequate semigroup with semilattice $E^0$ and let $I = \cup_{x\in E^0}{L_x}$ be a left regular band with a semilattice transversal $E^0$. Suppose also that there is defined on $I$ a left $S^0-$action $S^0\times I \to I$ given by $(x,e) \mapsto x\ast e$ and which is distributive over the multiplication on $I$ satisfying:
\begin{enumerate}
\item for all $x,y\in S^0$, $x\ast y^+ = (xy)^+$,

\item if $x,x_1, x_2\in S^0, e_1\in L_{x_1^+}, e_2\in L_{x_2^+}$ and if
$$x^+(x\ast e_1) = x^+(x\ast e_2),\ xx_1 = xx_2$$
then
$$x^\ast\ast e_1= x^\ast\ast e_2,\ x^\ast x_1 = x^\ast x_2$$
\end{enumerate}
Define a multiplication on the set
$$
W=\{(e,x)\in I\times S^0 : e \in L_{x^+}\}
$$
by
$$(e,x)(g,y) = (e(x\ast g),xy).
$$
Then $W$ is a left adequate, quasi-adequate semigroup with an admissible, left ample, adequate transversal isomorphic to $S^0$.

If in addition we have
\begin{enumerate}
\item[3.] $x^+\ast e = x^+e$ for all $e \in I,x\in S^0$
\end{enumerate}
then $I(W)\cong I$.

Moreover every left adequate, quasi-adequate semigroup $S$ with a left ample, admissible adequate transversal can be constructed in this way.
\end{theorem}

\begin{proof}
We show that the conditions for Theorem~\ref{structure-theorem} are satisfied. Let $I=\bigcup_{x\in E^0}{L_x}, E^0$ and $S^0$ be as in the statement of the theorem and let $\Lambda = E^0 = \bigcup_{x\in E^0}{R_x}$. Let $x,y \in S^0$ and define $\alpha_{x,y} : R_{x^\ast}\times L_{y^+}\to L_{(xy)^+}$ by $(x^\ast,e)\alpha_{x,y} = x\ast e$ and define $\beta_{x,y} : R_{x^\ast}\times L_{y^+}\to R_{(xy)^\ast}$ by $(x^\ast,e)\beta_{x,y} = (xy)^\ast$.

\begin{enumerate}
\item if $f \in R_{x^\ast}, g \in L_{y^+}, h \in R_{y^\ast}, k \in L_{z^+}$ then
$$
\begin{array}{rl}
(f,g)\alpha_{x,y}\left((f,g)\beta_{x,y}h,k\right)\alpha_{xy,z}&=(x\ast g)((xy)\ast k)\\
&=(x\ast g)(x\ast (y\ast k))\\
&=x\ast\left(g(y\ast k)\right)\\
&=\left(f,g(h,k)\alpha_{y,z}\right)\alpha_{x,yz}
\end{array}
$$
and
$$
\begin{array}{rl}
\left(f,g(h,k)\alpha_{y,z}\right)\beta_{x,yz}(h,k)\beta_{y,z}& = (x(yz))^\ast(yz)^\ast\\
&=((xy)z)^\ast\\
&=\left((f,g)\beta_{x,y}h,k\right)\beta_{xy,z}.
\end{array}
$$

\item $(x^\ast,y^+)\alpha_{x,y} = x\ast y^+ = (xy)^+$.
\item if $x, x_1, x_2\in S^0, e_1\in L_{x_1^+}, f_1 \in R_{x_1^\ast}, e_2\in L_{x_2^+}, f_2 \in R_{x_2^\ast}, e\in L_{x^+}$ and if
$$e_1(f_1,e)\alpha_{x_1,x} = e_2(f_2,e)\alpha_{x_2,x},\ x_1x = x_2x \text{ and } (f_1,e)\beta_{x_1,x}x^\ast = (f_2,e)\beta_{x_2,x}x^\ast$$
then $(f_1,e)\alpha_{x_1,x} = x_1\ast e = (f_1,e)\alpha_{x_1,x^+}$ and $x_1x^+ = x_2x^+$ as $x\rs x^+$. Hence
$$e_1(f_1,e)\alpha_{x_1,x^+} = e_2(f_2,e)\alpha_{x_2,x^+},\ x_1x^+ = x_2x^+ \text{ and } (f_1,e)\beta_{x_1,x^+} = (f_2,e)\beta_{x_2,x^+}.$$
\item if $x,x_1, x_2\in S^0, e_1\in L_{x_1^+}, f_1 \in R_{x_1^\ast}, e_2\in L_{x_2^+}, f_2 \in R_{x_2^\ast}, f\in R_{x^\ast}$ and if
$$x^+(f,e_1)\alpha_{x,x_1} = x^+(f,e_2)\alpha_{x,x_2},\ xx_1 = xx_2 \text{ and } (f,e_1)\beta_{x,x_1}f_1 = (f,e_2)\beta_{x,x_2}f_2$$
then $x^+(x\ast e_1) = x^+(x\ast e_2)$ and so $x^\ast\ast e_1 = x^\ast\ast e_2$ and $x^\ast x_1 = x^\ast x_2$ by property (2). Hence
$$(f,e_1)\alpha_{x^\ast,x_1} = (f,e_2)\alpha_{x^\ast,x_2},\ x^\ast x_1 = x^\ast x_2 \text{ and } (f,e_1)\beta_{x^\ast,x_1}f_1 = (f,e_2)\beta_{x^\ast,x_2}f_2.$$
\end{enumerate}

We can consequently deduce from Theorem~\ref{structure-theorem} that $W = \{(e,x,f)\in I\times S^0\times\Lambda : e\in L_{x^+}, f\in R_{x^\ast}\}$ is a quasi-adequate semigroup with an admissible adequate transversal $W^0=\{(x^+,x,x^\ast) : x\in S^0\}\cong S^0$. It is also clear that $W\cong\{(e,x) \in I\times S^0 : e \in L_{x^+}\}$.
Hence from Lemma~\ref{left-adequate-lemma}
we see that $W$ is left adequate and from Lemma~\ref{left-ample-lemma} that 
$W^0$ is left ample as required.

Suppose that $x^+\ast e = x^+e$ for all $e \in I,x\in S^0$ and consider the map $I\to W$ given by $a\mapsto(a,a^0)$, where $a^0$ is the unique inverse of $a$ in the semilattice transversal $E^0$. It is not too difficult to see that in fact $a=aa^0$ and $a^0=a^0a$ and so $(a^0)^+ =a^0\l x$. Now $ab\mapsto(ab,(ab)^0) = (a(a^0b),(ab)^0) = (a(a^0\ast b),a^0b^0) = (a,a^0)(b,b^0)$ and so $I(W)\cong I$.
\medskip

Conversely, let $S$ be a left adequate, quasi-adequate semigroup with a left ample, admissible, adequate transversal $S^0$.
\smallskip

From Theorem~\ref{structure-theorem} we see that $S\cong W$ where
$$
W=\{(e,x,f)\in I\times S^0\times\Lambda : e \in L_{x^+}, f \in R_{x^\ast}\}
$$
and multiplication is given by
$$(e,x,f)(g,y,h) = (e(f,g)\alpha_{x,y},xy,(f,g)\beta_{x,y}h) = (ee_{xfgy},xy,f_{xfgy}h).
$$
\smallskip
Since $S$ is left adequate we in fact have
$$
W\cong\{(e,x)\in I\times S^0 : e \in L_{x^+}\}
$$
and multiplication is given by
$$(e,x)(g,y) = (ee_{xx^\ast gy},xy) = (ee_{xgy},xy).
$$
Let $x\in S^0$ and define a left action of $S^0$ on $I$ by $x\ast e = e_{xey}$ where $e \in L_{y^+}$. To check that $\ast$ is indeed an action notice first that $e_{xey}\rs xey\rs xey^+ = xe\rs e_{xe}$ and so $e_{xey} = e_{xe}$ and $x\ast e$ is independent of $y$. Hence $\ast$ is well defined. Now let $x,y \in S^0, e \in I$.
Then $e_{ye}\rs ye$ and so $xe_{ye}\rs xye$ and hence from Lemma~\ref{rs-ls-lemma} we see that $x\ast\left(y\ast e\right)= x\ast\left(e_{ye}\right) = e_{xe_{ye}} = e_{xye} = (xy)\ast e$. Consequently $\ast$ is a left action.

\smallskip

Suppose $e\in L_{y^+}, f\in L_{z^+}$ so that $ef\in L_{y^+}L_{z^+} = L_{y^+z^+}=L_{(z^+y^+)^+}=L_{(z^+y)^+}$. Hence $x\ast(ef) = e_{xef}$. Now
$(x\ast e)(x\ast f) = e_{xe}e_{xf} = e_{e_{xe}e_{xf}}$ and since $e_{xf}\rs xf$ then $e_{xe}e_{xf}\rs e_{xe}{xf}$ and so from Lemma~\ref{rs-ls-lemma} we see that $e_{xe}e_{xf} = e_{e_{xe}e_{xf}} = e_{e_{xe}{xf}}$.
But
$$
\begin{array}{lcl}
xe& = &e_{xe}xe = e_{xe}e_{\overline{xe}}xe = e_{xe}e_{\overline{xe}}xef_{\overline{xe}}\\
&=&e_{xe}\overline{xe} = e_{xe}\x\;\e = e_{xe}\left(\x\;\e\right)^+\x\\
&=&e_{xe}e_{\overline{xe}}\x = e_{xe}x.\\
\end{array}
$$

Hence $e_{xef} = e_{e_{xe}xf} = e_{xe}e_{xf}$ and so $\ast$ satisfies the appropriate distributive property.

\medskip

Notice that if $x\in S^0$ then $e_x\rs x\rs x^+$ and so $e_x = x^+$ since $S$ is left adequate. Let $x,y \in S^0$. Then $x\ast y^+ = e_{xy^+} = \left(xy^+\right)^+ = (xy)^+$ and so property (1) holds.

\medskip

Suppose that $x,x_1, x_2\in S^0, e_1\in L_{x_1^+}, e_2\in L_{x_2^+}$ and that
$$x^+(x\ast e_1) = x^+(x\ast e_2),\ xx_1 = xx_2.$$
Then $x^+e_{xe_1} = x^+e_{xe_2},\ xx_1 = xx_2$ and since $x\ls x^\ast$ then $x^\ast x_1 = x^\ast x_2$.

Now $e_{xe_1}\rs xe_1$ and so $x^+e_{xe_1}\rs x^+xe_1 = xe_1\rs e_{xe_1}$. Hence $x^+e_{xe_1} = e_{xe_1}$ since $|I\cap R^\ast_{xe_1}|=1$. Consequently we deduce that $e_{xe_1} = e_{xe_2}$. Now using Theorem~\ref{quasi-adequate-ef-theorem} we deduce that $e_{xe_1x_1} = e_{xe_1}e_{\overline{xe_1}e_{x_1}x_1} = e_{xe_1}e_{xx_1^+x_1} = e_{xe_1}e_{xx_1} = e_{xe_2}e_{xx_2} = e_{xe_2x_2}$. Hence $xe_1x_1 = e_{xe_1x_1}\overline{xe_1x_1} = e_{xe_1x_1}xx_1 = xe_2x_2$. Since $x\ls x^\ast$ then $x^\ast e_1x_1 = x^\ast e_2x_2$. Consequently $e_{x^\ast e_1x_1}=e_{x^\ast e_2x_2}$ and from above we see that $x\ast e_1 = e_{x^\ast e_1}=e_{x^\ast e_1x_1}=e_{x^\ast e_2x_2}=e_{x^\ast e_2} = x\ast e_2$ and so property (2) holds.

Finally if $x \in S^0,e \in I$ then $x^+\ast e = e_{x^+e} = x^+e$ and so property (3) holds.
\end{proof}

\section{Inverse Transversals}

We now consider the situation when $S$ is in fact regular. We make use of Theorem~\ref{main-regular-theorem} together with the associated remarks.
\begin{proposition}
Let $S$ be a quasi-adequate semigroup with an adequate transversal $S^0$. Then $S$ is orthodox if and only if $S^0$ is inverse.
\end{proposition}

\begin{proof}
If $S$ is orthodox then clearly $S^0$ is regular and so inverse. Conversely if $S^0$ is inverse then for all $x\in S$ we have $x = e_x\x f_x \in \langle Reg(S)\rangle$ and so from~\cite[Proposition 1.3]{fountain3} $x \in Reg(S)$.
\end{proof}

Let $S^0,I,\Lambda$ and $W$ be as in the statement of Theorem~\ref{structure-theorem} together with the given multiplication and suppose that properties (1) and (2) hold so that $W$ is a semigroup.
Let $(e,x,f) \in W$ and notice that $e\in L_{x^+}, f\in R_{x^\ast}$ and so $e^{00} = \overline{e} = x^+, f^{00} = \overline{f} = x^\ast$.
Now suppose in what follows that $S^0$ is an inverse semigroup so that on $S^0$, ${\cal R}^\ast = {\cal R}$. Notice that $x^+ = xx^{-1}$ and $x^\ast = x^{-1}x$. Also it is easy to see that

\begin{enumerate}
\item $\left(x^{-1}\right)^+ = x^{-1}\left(x^{-1}\right)^{-1} = x^{-1}x = x^\ast$,
\item $\left(x^{-1}\right)^\ast = \left(x^{-1}\right)^{-1}x^{-1} = xx^{-1} = x^+$,
\item $\left(x^+\right)^{-1}= x^+$,
\item$\left(x^\ast\right)^{-1}= x^\ast$.
\end{enumerate}

From above and from Lemma~\ref{E0-transversal-lemma} we see that $e^0 = e^{00} = x^+=(x^{-1})^\ast$ and similarly $f^0 = f^{00} = x^\ast=(x^{-1})^+$
and consequently $(f^0,x^{-1},e^0) \in W$. It is then easy to check that
$$
(e,x,f)(f^0,x^{-1},e^0) = (e,x^+,x^+)=(e,e^0,e^0)
$$
and that
$$
(e,e^0,e^0)(e,x,f) = (e,x,f).
$$
Hence it follows that $(e,e^0,e^0)\r(e,x,f)$ and that $W$ is regular. It is also easy to verify that $(f^0,x^{-1},e^0)\in V(e,x,f)$.

Now suppose that the condition in property (3) of Theorem~\ref{structure-theorem} holds, namely
$$
e_1(f_1,e)\alpha_{x_1,x} = e_2(f_2,e)\alpha_{x_2,x},\ x_1x = x_2x \text{ and } (f_1,e)\beta_{x_1,x}x^\ast = (f_2,e)\beta_{x_2,x}x^\ast.
$$
Then we have
$$
(e_1,x_1,f_1)(e,x,x^\ast)=(e_2,x_2,f_2)(e,x,x^\ast)
$$
and so multiplying on the right by $(\left(x^\ast\right)^0,x^{-1},e^0)$ we see that
$$
(e_1,x_1,f_1)(e,x^+,x^+)=(e_2,x_2,f_2)(e,x^+,x^+).
$$
Consequently we see that the conclusion of property (3) holds. In a similar way property (4) holds as well.
Hence we have proved

\begin{corollary}\label{orthodox-structure-theorem}{\ \rm (\cite[Theorem 3.6]{saito})}
Let $S^0$ be an inverse semigroup with semilattice of idempotents $E^0$ and let $I = \cup_{x\in E^0}{L_x}$ be a left regular band and $\Lambda = \cup_{x\in E^0}{R_x}$ a right regular band with a common semilattice transversal $E^0$. Suppose that for each $x,y \in S^0$ there exist $\alpha_{x,y}\in PT(\Lambda\times I,I)$ and $\beta_{x,y}\in PT(\Lambda\times I,\Lambda)$ satisfying:
\begin{enumerate}
\item $\text{\rm dom}(\alpha_{x,y}) = \text{\rm dom}(\beta_{x,y}) \subseteq R_{x^{-1}x}\times L_{yy^{-1}}, (f,e)\alpha_{x,y} \in L_{(xy)(xy)^{-1}}$ and $(f,e)\beta_{x,y}\in R_{(xy)^{-1}(xy)}$,
\item if $f \in R_{x^{-1}x}, g \in L_{yy^{-1}}, h \in R_{y^{-1}y}, k \in L_{zz^{-1}}$ then
$$(f,g)\alpha_{x,y}\left((f,g)\beta_{x,y}h,k\right)\alpha_{xy,z} = \left(f,g(h,k)\alpha_{y,z}\right)\alpha_{x,yz}$$
$$\left(f,g(h,k)\alpha_{y,z}\right)\beta_{x,yz}(h,k)\beta_{y,z} = \left((f,g)\beta_{x,y}h,k\right)\beta_{xy,z},$$
\item $(x^{-1}x,yy^{-1})\alpha_{x,y} = (xy)(xy)^{-1}, (x^{-1}x,yy^{-1})\beta_{x,y} = (xy)^{-1}(xy)$,
\end{enumerate}
Define a multiplication on the set
$$
W=\{(e,x,f)\in I\times S^0\times\Lambda : e \in L_{xx^{-1}}, f \in R_{x^{-1}x}\}
$$
by
$$(e,x,f)(g,y,h) = (e(f,g)\alpha_{x,y},xy,(f,g)\beta_{x,y}h).
$$
Then $W$ is an orthodox semigroup with an inverse transversal isomorphic to $S^0$. Moreover, if in addition $\alpha$ and $\beta$ satisfy
\begin{enumerate}
\item[4.] for all $f\in\Lambda, e \in I$,
$$(f^0,e)\alpha_{f^0,e^0} = f^0e,\quad (f,e^0)\beta_{f^0,e^0} = fe^0,$$
\end{enumerate}
then $I(W) \cong I, \Lambda(W)\cong \Lambda$.

Conversely every orthodox semigroup $S$, with an inverse transversal can be constructed in this way.
\end{corollary}

From Theorem~\ref{left-adequate-theorem} we can also deduce the following

\begin{corollary}{\ \hbox{\rm (Cf. \cite[Theorem 1]{yoshida})}}
Let $S^0$ be an inverse semigroup with semilattice of idempotents $E^0$ and let $I$ be a left regular band with a semilattice transversal isomorphic to $E^0$. Suppose that on $I$ we have a left action of $S^0$ given by $(x,e)\mapsto x\ast e$ and which is distributive over the multiplication on $I$ satisfying
\begin{enumerate}
\item for all $x,y \in S^0, x\ast (yy^{-1}) = (xy)(xy)^{-1}$;
\item for all $x \in S^0, e \in I, (xx^{-1})\ast e = (xx^{-1})e$.
\end{enumerate}
Define a multiplication on the set
$$
W=\{(e,x)\in I\times S^0 : e \in L_{xx^{-1}}\}
$$
by
$$(e,x)(g,y) = (e(x\ast g),xy).
$$
Then $W$ is a left inverse semigroup with an inverse transversal isomorphic to $S^0$ and $I(W)\cong I$.
Moreover every left inverse semigroup $S$ with an inverse transversal can be constructed in this way.
\end{corollary}

\begin{proof}
Let $I, S^0,E^0$ and $\ast$ be as given and suppose that $x,x_1,x_2 \in S^0, e_1 \in L_{x_1x_1^{-1}}, e_2 \in L_{x_2x_2^{-1}}$ and that
$$
xx^{-1}(x\ast e_1) = xx^{-1}(x\ast e_2), xx_1 = xx_2.
$$
Then $(xx^{-1},x)(e_1,x_1) = (xx^{-1},x)(e_2,x_2)$ and so multiplying on the left by $((x^\ast)^0,x^{-1})$ we see that
$(x^\ast\ast e_1,x^\ast x_1)=(x^\ast\ast e_2,x^\ast x_2)$. The result now follows now by Theorem~\ref{left-adequate-theorem}.
\end{proof}

From Corollary~\ref{quasi-ideal-corollary} we can deduce

\begin{corollary}
Let $S^0$ be an inverse semigroup with semilattice $E^0$ and let $I = \cup_{x\in E^0}{L_x}$ be a left normal band and $\Lambda = \cup_{x\in E^0}{R_x}$ a right normal band with a common semilattice transversal $E^0$. Let
$$
W=\{(e,x,f)\in I\times S^0\times\Lambda : e \in L_{xx^{-1}}, f \in R_{x^{-1}x}\}
$$
and define a multiplication by
$$
(e,x,f)(g,y,h) = (e(xy)(xy)^{-1},xy,(xy)^{-1}(xy) h).
$$
Then $W$ is an orthodox semigroup with a quasi-ideal inverse transversal $T^0\cong S^0$. Conversely every such transversal can be constructed in this way.
\end{corollary}

Finally from Corollary~\ref{spined-product-corollary} we deduce
\begin{corollary}
Let $L$ be a left inverse semigroup and $R$ a right inverse semigroup with a common quasi-ideal inverse transversal $S^0$. Construct the spined product
$$
\spd = \{(x,a)\in L\times R : x^0 = a^0\}
$$
and define a multiplication on $\spd$ by
$$
(x,a)(y,b) = (xy^{00},a^{00}b)
$$
Then $\spd$ is an orthodox semigroup with a quasi-ideal inverse transversal isomorphic to $S^0$. Moreover every such transversal can be constructed in this way.
\end{corollary}

The authors would like to thank V. Gould and J. Fountain for useful discussions relating to this work.

\end{document}